\documentclass[3p,10pt,a4paper,twoside,fleqn,sort&compress]{arxiv}

\usepackage{amssymb,amsmath,latexsym}
\usepackage[varg]{pxfonts}

\newtheorem{theorem}{Theorem}[section]

\newtheorem{lemma}[theorem]{Lemma}

\newcommand{\FilPages}{}

\newcommand{\FilReceived}{}

\begin{document}

\title{A Note on the Reduction Formulas for \\ Some Systems of Linear Operator Equations}

\author[]{Ivana Jovovi\' c}
\ead{ivana@etf.rs}
\author[]{Branko Male\v sevi\' c}
\ead{malesevic@etf.rs}
\address{Faculty of Electrical Engineering, University of Belgrade}
\newcommand{\AuthorNames}{I. Jovovi\' c, B. Male\v sevi\' c}
\newcommand{\FilMSC}{15A21}
\newcommand{\FilKeywords}{Partial and total reduction of linear systems of operator equations,
the rational and the Jordan canonical forms, doubly companion matrix,
characteristic polynomial, characteristic matrix, adjugate matrix}
\newcommand{\FilCommunicated}{Professor Dragana Cvetkovi\' c Ili\' c}

\begin{abstract}
We consider partial and total reduction of a nonhomogeneous linear system of the operator equations
with the system matrix in the same particular form as in paper \cite{ShayanfarHadizadeh2013}.
Here we present two different concepts. One is concerned with partially reduced systems
obtained by using the Jordan and the rational form of the system matrix.
The other one is dealing with totally reduced systems obtained by finding
the adjugate matrix of the characteristic matrix of the system matrix.
\end{abstract}

\maketitle
\makeatletter
\renewcommand\@makefnmark%
{\mbox{\textsuperscript{\normalfont\@thefnmark)}}}
\makeatother

\section{Introduction}

Let $K$ be a field, $V$ a vector space over $K$ and $A : V \rightarrow V$ a linear operator on the vector space $V$.
Nonhomogeneous linear system of the operator equations with constant coefficients $b_{ij} \in K$ in unknowns $x_i$, $1 \leq i,j \leq n$,
has the general form
\begin{equation}
\label{generalsystem}
\begin{array}{lcl}
A(x_{1})   & = & b_{11}x_1 + b_{12}x_2 + \ldots + b_{1n}x_n + \varphi_1       \\[0.5 ex]
A(x_{2})   & = & b_{21}x_1 + b_{22}x_2 + \ldots + b_{2n}x_n + \varphi_2       \\
\;\vdots   &   &                                                              \\
A(x_{n})   & = & b_{n1}x_1 + b_{n2}x_2 + \ldots + b_{nn}x_n + \varphi_n,
\end{array}
\end{equation}
for $\varphi_1, \varphi_2, \ldots , \varphi_n \in V$.
It is common to consider system (\ref{generalsystem}) in the matrix form
\[
\vec{A}(\vec{x})= B\vec{x} + \vec{\varphi},
\]
where $\vec{x}=[x_{1} \,  x_{2} \, \ldots \, x_{n}]^T \in V^{n \times 1}$ is a column of unknowns,
$\vec{A}: V^{n \times 1} \rightarrow V^{n \times 1}$ is a vector operator defined componentwise by
$\vec{A}(\vec{x})=[A(x_1) \, A(x_2) \, \ldots \, A(x_n)]^T$,
$\vec{\varphi}=[\varphi_{1} \, \varphi_{2} \, \ldots  \,\varphi_{n}]^T \in V^{n \times 1}$ is a nonhomogeneous term and
$B=[b_{ij}]_{i,j=1}^n \in K^{n \times n}$ is a system matrix.
There is an extensive amount of literature relating to a problem of reducing the linear system of the operator equations
to an equivalent system in a simpler form.
The most widely used technique is to transform the system into a system with block diagonal system matrix using basis transformation.
Such a basis transformation produces a new system equivalent to the initial one, which is decoupled into several subsystems.
Each subsystem corresponds to one block of the new system matrix.
Subsystems are uncoupled, so we may solve them separately,
and then simply assemble these individual solutions together to obtain a solution of the general system.
The Jordan canonical form is the most commonly used if the field of coefficients is algebraically closed.
If it is not the case, the Jordan canonical form of a matrix can only be achieved by adding a field extension.
The rational canonical form of a matrix is the best diagonal block form that can be achieved over the field of coefficients
and it corresponds to the factorization of the characteristic polynomial into invariant factors without adding any field extension.
In the paper \cite{MalesevicTodoricJovovicTelebakovic2010} the idea has been to use the rational canonical form
to reduce the linear system (\ref{generalsystem}) to an equivalent partially reduced system.
Each subsystem of the partially reduced system consists of a higher order linear operator equation having only one variable
and the first order linear operator equations in two variables.
By the order of a linear operator equation we mean the highest power of the operator in the equation.
An another method for solving the linear system of the operator equations,
which does not require a change of basis, is discussed in \cite{MalesevicTodoricJovovicTelebakovic2012}.
The system is reduced to a so called totally reduced system, i.e. to a system with separated variables,
by using the characteristic polynomial $\Delta_{B}(\lambda)= \det(\lambda I - B)$ of the system matrix $B$.
This system consists of a higher order operator equations which differ only in the variables and in the nonhomogeneous terms.
The general reduction formulas from \cite{MalesevicTodoricJovovicTelebakovic2010} and \cite{MalesevicTodoricJovovicTelebakovic2012}
can be applied to some special systems.
An example of total reduction of the linear system of the operator equations with the system matrix in the companion form
can be found in \cite{Jovovic2013}.

In this paper we consider nonhomogeneous linear systems of the operator equations in the form of Shayanfar and Hadizadeh \cite{ShayanfarHadizadeh2013}:
\begin{subequations}
\begin{equation}
\tag{$\hat{2}$}
\label{system}
\begin{array}{lcl}
A(x_{1})   & = & b_1x_1 + b_1x_2 + \ldots + b_1x_n + \varphi_1       \\[0.5 ex]
A(x_{2})   & = & b_2x_1 + b_2x_2 + \ldots + b_2x_n + \varphi_2       \\
\;\vdots   &   &
\\
A(x_{n})   & = & b_nx_1 + b_nx_2 + \ldots + b_nx_n + \varphi_n
\end{array}
\end{equation}
\;\;\mbox{and}\;\;
\begin{equation}
\tag{$\check{2}$}
\label{transpose system}
\begin{array}{lcl}
A(x_{1})   & = & b_1x_1 + b_2x_2 + \ldots + b_nx_n + \varphi_1       \\[0.5 ex]
A(x_{2})   & = & b_1x_1 + b_2x_2 + \ldots + b_nx_n + \varphi_2       \\
\;\vdots   &   &
\\
A(x_{n})   & = & b_1x_1 + b_2x_2 + \ldots + b_nx_n + \varphi_n,
\end{array}
\end{equation}
\end{subequations}
for $b_1, b_2, \ldots b_n \in K$ and $\varphi_1, \varphi_2, \ldots , \varphi_n \in V$.
The systems (\ref{system}) and (\ref{transpose system}) can be rewritten in the matrix form
$$
\vec{A}(\vec{x})=\hat{\mathstrut B}\vec{x} + \vec{\varphi} \;\;\mbox{and}\;\; \vec{A}(\vec{x})=\check{\mathstrut B}\vec{x} + \vec{\varphi},
$$
where
$$
\hat{\mathstrut B}=
\left[
\begin{array}{cccc}
b_1    & b_1    & \ldots & b_1    \\
b_2    & b_2    & \ldots & b_2    \\
\vdots & \vdots & \ddots & \vdots \\
b_n    & b_n    & \ldots & b_n
\end{array}
\right]
\;\mbox{and}\;
\check{\mathstrut B}=
\left[
\begin{array}{cccc}
b_1    & b_2    & \ldots & b_n    \\
b_1    & b_2    & \ldots & b_n    \\
\vdots & \vdots & \ddots & \vdots \\
b_1    & b_2    & \ldots & b_n
\end{array}
\right].
$$
Shayanfar and Hadizadeh in \cite{ShayanfarHadizadeh2013} have used matrix polynomial approach to reduce the systems (\ref{system}) and (\ref{transpose system}).
Their main idea has been to apply the Smith canonical form for obtaining a class of independent equations.
Here we present how standard Jordan and rational canonical forms can help in solving this problem.
We also sketch a method for total reduction for these special cases.

\newpage

$\;$

\section{Some Basic Notions and Notations}

In this section we will review some standard facts from linear algebra, see for example \cite{DummitFoote04, Gantmacher00}.

Let $B$ be an $n \times n$ matrix over the field $K$. An element $\lambda \in K$ is called an eigenvalue of $B$ with corresponding eigenvector $v$,
if $v$ is a nonzero $n \times 1$ column with entries in $K$ such that $\lambda v = B v$.
The set of all eigenvectors with the same eigenvalue $\lambda$, together with the zero vector,
is a vector space called the eigenspace of the matrix $B$ that corresponds to the eigenvalue $\lambda$.
The geometric multiplicity of an eigenvalue $\lambda$ is defined as the dimension of the associated eigenspace, i.e.
it is the number of linearly independent eigenvectors corresponding to that eigenvalue.
The algebraic multiplicity of an eigenvalue $\lambda$ is defined as the multiplicity
of the corresponding root of the characteristic polynomial $\Delta_{B}(\lambda)= \det(\lambda I - B)$.
A generalized eigenvector $u$ of $B$ associated to $\lambda$ is a nonzero $n \times 1$ column with entries in $K$ satisfying $(B - \lambda I)^k u = 0$,
for some $k \in \mathbb{N}$.
The set of all generalized eigenvectors for a given eigenvalue $\lambda$,
together with the zero vector, form the generalized eigenspace for $\lambda$.

Recall that the $n \times n$ matrix $\widetilde{B}$ over the field $K$ is said to be similar to the matrix $B$
if there exists an $n \times n$ nonsingular matrix $S$ with entries in $K$ such that $\widetilde{B} = S^{-1} B S$.
The matrix $S$ is called transition matrix.
Similarity of matrices is an equivalence relation on $K^{n \times n}$.
It can be shown that basic notions of linear algebra such as rank, characteristic polynomial, determinant, trace,
eigenvalues and their algebraic and geometric multiplicities, minimal polynomial,
the Jordan and the rational canonical forms, elementary divisors and invariant factors are similarity invariants.

One can notice that $\check{\mathstrut B}=\hat{\mathstrut B}^T$.
Since a matrix and its transpose matrix are similar,
and all characteristic that we examine are similarity invariant,
we focus only on the first system (\ref{system}).
All conclusions we make are the same for the second one as well.

Let us consider more closely some similarity invariants.
Invariant factors of the matrix $B$ are polynomials
$$\Delta_{1}(\lambda)=\frac{D_{1}(\lambda)}{D_{0}(\lambda)}, \Delta_{2}(\lambda)=\frac{D_{2}(\lambda)}{D_{1}(\lambda)}, \ldots ,\Delta_{r}(\lambda)=\frac{D_{r}(\lambda)}{D_{r-1}(\lambda)},$$
where $D_{j}(\lambda)$ is the greatest common divisor of all minors of the order $j$ in the matrix $\lambda I - B$ and $D_{0}(\lambda)=1$, $1 \leq j \leq r$.
Elementary divisors of the matrix $B$ are monic irreducible polynomials over the field $K$ into which the invariant factors of the matrix $B$ split.
The product of all elementary divisors of the matrix $B$ is its characteristic polynomial, and their least common multiple is its minimum polynomial.

The $k \times k$ matrix of the form
$$J \!=\! \left[
\begin{array}{ccccc}
\lambda  & 1        & \ldots & 0       & 0          \\
0        & \lambda  & \ldots & 0       & 0          \\
\vdots   & \vdots   &        & \vdots  & \vdots     \\
0        & 0        & \ldots & \lambda & 1          \\
0        & 0        & \ldots & 0       & \lambda
\end{array}
\right]$$
is called the Jordan block of size $k$ with eigenvalue $\lambda$.
A matrix is said to be in the Jordan canonical form if it is a block diagonal matrix with Jordan blocks along the diagonal.
The number of Jordan blocks corresponding to an eigenvalue $\lambda$ is equal to its geometric multiplicity
and the sum of their sizes is equal to the algebraic multiplicity of $\lambda$.
Every square matrix is similar to a matrix in the Jordan canonical form.
That is, for the matrix $B$ there exists a nonsingular matrix $S$,
so that $J = S^{-1} B S$, where $J$ is in the Jordan canonical form.
The Jordan canonical form of the matrix $B$ is unique up to the order of the Jordan blocks.

If the matrix $B$ is similar to a diagonal matrix, then $B$ is said to be diagonalizable.
For the diagonalizable matrix $B$, $J = S^{-1} B S$ is diagonal matrix,
$S$ is a matrix obtained by using as its columns any set of linearly independent eigenvectors of $B$,
and the diagonal entries of $J$ are their corresponding eigenvalues.
The matrix $S$ is called modal matrix.

The companion matrix of the polynomial $\Delta(\lambda)= \lambda^n + d_{1} \lambda^{n-1} + \ldots + d_{n-1} \lambda + d_n$ is the matrix
$$C \!=\! \left[
\begin{array}{cccccc}
0        & 1       & 0       & \ldots & 0       & 0       \\
0        & 0       & 1       & \ldots & 0       & 0       \\
\vdots   & \vdots  & \vdots  & \ddots & \vdots  & \vdots  \\
0        & 0       & 0       & \ldots & 0       & 1       \\
-d_{n}   &-d_{n-1} &-d_{n-2} & \ldots &-d_{2}   &-d_{1}
\end{array}
\right].$$
It can easily be seen that the characteristic polynomial of the companion matrix $C$ is $\Delta(\lambda)$.
A matrix is said to be in the rational canonical form if it is a block diagonal matrix with blocks the companion matrices of the monic polynomials
$\Delta_1(\lambda), \Delta_2(\lambda), \ldots , \Delta_r(\lambda)$ of degree at least one with $\Delta_1(\lambda) \; | \; \Delta_2(\lambda) \; | \; \ldots \; | \; \Delta_r(\lambda)$.
Every square matrix is similar to a matrix in the rational canonical form.
That is, for the matrix $B$ there exists a nonsingular matrix $T$,
so that $C = T^{-1} B T$, where $C$ is in the rational canonical form.
The rational canonical form of the matrix $B$ is unique up to the order of the diagonal blocks.

Butcher and Chartier in \cite{ButcherWright2006} introduced the notion of the doubly companion matrix of polynomials
$\alpha(\lambda)=\lambda^n + a_{1}\lambda^{n-1} + \ldots + a_{n-1}\lambda + a_{n}$ and
$\beta(\lambda)=\lambda^n + b_{1}\lambda^{n-1} + \ldots + b_{n-1}\lambda + b_{n}$
as an $n \times n$ matrix over the field $K$ of the form
$$
C(\alpha, \beta)=
\left[
\begin{array}{ccccc}
-a_{1}      & -a_{2}      & \ldots       & -a_{n-1}      & -a_{n}-b_{n}       \\
1           & 0           & \ldots       & 0             & -b_{n-1}           \\
\vdots      & \vdots      & \ddots       & \vdots        & \vdots             \\
0           & 0           & \ldots       & 0             & -b_{2}             \\
0           & 0           & \ldots       & 1             & -b_{1}
\end{array}
\right].
$$

\section{Preliminaries and Auxiliary Results}

At the beginning of this section we first present two standard lemmas for reduction process using canonical forms.
Afterwards we give a brief exposition of two main results from \cite{MalesevicTodoricJovovicTelebakovic2010} and \cite{MalesevicTodoricJovovicTelebakovic2012}.

\begin{lemma}
\label{Jordan form}
Let $J$ be the Jordan canonical form of the matrix $B$,
i.e. there exists a nonsigular matrix $S$ such that $J=S^{-1} B S$.
Then the system {\rm (\ref{generalsystem})} given in the matrix form
$$
\vec{A}(\vec{x})=B\vec{x} + \vec{\varphi},
$$
can be reduced to the system
$$
\vec{A}(\vec{y})=J\vec{y} + \vec{\psi},
$$
where $\vec{\psi}=S^{-1}\vec{\varphi}$ is its nonhomogeneous term and $\vec{y}=S^{-1}\vec{x}$ is a column of the unknowns.
\end{lemma}

\begin{lemma}
\label{rational form}
Let $C$ be the rational canonical form of the matrix $B$, i.e.
there exists a nonsigular matrix $T$ such that $C=T^{-1}BT$.
Then the system {\rm (\ref{generalsystem})} given in the matrix form
$$
\vec{A}(\vec{x})=B\vec{x} + \vec{\varphi},
$$
can be reduced to the system
$$
\vec{A}(\vec{z})=C\vec{z} + \vec{\nu},
$$
where $\vec{\nu}=T^{-1}\vec{\varphi}$ is its nonhomogeneous term and $\vec{z}=T^{-1}\vec{x}$ is a column of the unknowns.
\end{lemma}

Let $\delta_{k}^{1}{\big (} B; \vec{A}^{n-k}(\vec{\varphi}){\big )}$ stands for the sum of the principal minors of the order $k$
containing the entries of the first column of the matrix obtained from the matrix $B$ by replacing column $\vec{A}^{n-k}(\vec{\varphi})$
in the place of the first column of $B$.
Following two theorems are concerned with partial reduction of the system (\ref{generalsystem}).

\begin{theorem}$($Theorem 3.4 from {\rm \cite{MalesevicTodoricJovovicTelebakovic2010}}$)$
\label{Theorem PRS special case}
Let us assume that the rational canonical form of the system matrix $B$ has only one block, i.e.
that the rational canonical form of $B$ is the companion matrix of the characteristic polynomial
$\Delta_{B}(\lambda) = \lambda^{n} + d_{1} \lambda^{n-1} + \ldots + d_{n-1} \lambda + d_{n}$.
Then the linear system of the operator equations {\rm (\ref{generalsystem})} can be transformed into equivalent partially reduced system
\[
\begin{array}{lcl}
\Delta_B(A)(y_{1})     & = & \sum_{k=1}^n (-1)^{k-1} \delta_{k}^{1}{\big (}C; \vec{A}^{n-k}(\vec{\psi}){\big )}     \\[0.5 ex]
y_{2}                  & = &  A(y_{1}) - \psi_{1}                                                                   \\[0.5 ex]
y_{3}                  & = &  A(y_{2}) - \psi_{2}                                                                   \\
\;\vdots               &   &                                                                                        \\
y_{n}                  & = &  A(y_{n-1}) - \psi_{n-1},
\end{array}
\]
where the columns $\vec{y}=[y_{1} \, y_{2} \, \ldots \, y_{n}]^{T}$ and $\vec{\psi}=[\psi_{1} \, \psi_{2} \, \ldots \, \psi_{n}]^{T}$
are determined by $\vec{y}=T^{-1}\vec{x}$ i $\vec{\psi}=T^{-1}\vec{\varphi}$ for a nonsingular matrix $T$ such that $C=T^{-1} B T$.
\end{theorem}
\begin{theorem}$($Theorem 3.7 from {\rm \cite{MalesevicTodoricJovovicTelebakovic2010}}$)$
\label{Theorem PRS}
Let us assume that the rational canonical form of the system matrix $B$ is block diagonal matrix
$$C=
\left[
\begin{array}{cccc}
C_{1}       & 0           & \ldots       & 0          \\
0           & C_{2}       & \ldots       & 0          \\
\vdots      & \vdots      & \ddots       & \vdots     \\
0           & 0           & \ldots       & C_{r}
\end{array}
\right]
\qquad (2 \leq k \leq n),$$
where $C_{i}$ are companion matrices of the monic polynomials
$\Delta_{C_{i}}(\lambda)=\lambda^{n_{i}} + d_{i,1}\lambda^{n_{i}-1} + \ldots + d_{i,n_{i}-1} \lambda + d_{i,n_{i}}$
of degree at least one with $\Delta_{C_{i}} \;|\; \Delta_{C_{i+1}}$ for $1 \leq i <r$.
Let $\ell_{1} \!=\! 0$ and $\ell_{i} \!=\! \sum_{j=1}^{i-1}{n_{j}}$ for $2 \leq\! i\! \leq k$.
Then the linear system of the operator equations {\rm (\ref{generalsystem})} can be transformed into equivalent partially reduced system
$\bigwedge\limits_{i=1}^{k}{(\mbox{$\cal R$}_{C_{i}})},$
where every subsystem $(\mbox{$\cal R$}_{C_{i}})$ is of the form
\[
\begin{array}{lcl}
\Delta_{C_{i}}(A)(y_{\ell_{i}+1}) & = & \sum_{k=1}^{n_{i}} (-1)^{k-1}\delta_{k}^{1}{\big (}C_{i};
[A^{n_i-k}(\psi_{\ell_{i}+1}) \, A^{n_i-k}(\psi_{\ell_{i}+2}) \, \ldots \, A^{n_i-k}(\psi_{\ell_{i}+n_{i}})]^{T}){\big )}                                                \\[0.5 ex]
y_{\ell_{i}+2}                    & = &  A(y_{\ell_{i}+1}) - \psi_{\ell_{i}+1}                                                                                           \\[0.5 ex]
y_{\ell_{i}+3}                    & = &  A(y_{\ell_{i}+2}) - \psi_{\ell_{i}+2}                                                                                           \\
\;\vdots                          &   &                                                                                                                                  \\
y_{\ell_{i}+n_{i}}                & = &  A(y_{\ell_{i}+n_{i}-1}) - \psi_{\ell_{i}+n_{i}-1}.
\end{array}
\]
The columns $\vec{y}=[y_{1} \, y_{2} \, \ldots \, y_{n}]^{T}$ and $\vec{\psi}=[\psi_{1} \, \psi_{2} \, \ldots \, \psi_{n}]^{T}$
are determined by $\vec{y}=T^{-1}\vec{x}$ i $\vec{\psi}=T^{-1}\vec{\varphi}$ for a nonsingular matrix $T$ such that $C=T^{-1} B T$.
\end{theorem}

We now present two theorems concerning total reduction of the system (\ref{generalsystem}).

\begin{theorem}$($Theorem 4.1 from {\rm \cite{MalesevicTodoricJovovicTelebakovic2010}}$)$
\label{Theorem adjugate matrix}
Assume that the system {\rm (\ref{generalsystem})} is given in the matrix form
$$
\vec{A}(\vec{x})=B\vec{x} + \vec{\varphi},
$$
and that matrices $B_{0}, \ldots, B_{n-1}$ are coefficients of the matrix polynomial $\lambda I - B$.
Then for the linear operator $\Delta_B (\vec{A})$, obtained by replacing $\lambda$ by $\vec{A}$
in the characteristic polynomial $\Delta_B(\lambda)$ of the system matrix $B$, the equality
$$
{\big (}\Delta_B (\vec{A}){\big )}(\vec{x}) = \sum_{k=1}^{n} B_{k-1}\vec{A}^{n-k}(\vec{\varphi})
$$
holds.
\end{theorem}

An explicit formula for total reduction of the system (\ref{generalsystem}) is given in the following theorem.

\begin{theorem} $($Theorem 4.3 from {\rm \cite{MalesevicTodoricJovovicTelebakovic2012}}$)$
\label{Theorem TRS}
Let $\delta_{k}^{i}{\big (} B; \vec{A}^{n-k}(\vec{\varphi}){\big )}$ stands for the sum of the principal minors of the order $k$
containing the entries of the $i^{th}$ column of the matrix obtained from the matrix $B$ by replacing column $\vec{A}^{n-k}(\vec{\varphi})$
in the place of the $i^{th}$ column of $B$.
Then the linear system of the operator equations {\rm (\ref{generalsystem})} implies the system,
which consists of the higher order operator equations as follows
\[
\begin{array}{lcl}
\Delta_B(A)(x_{1}) & = & \sum_{k=1}^n (-1)^{k-1} \delta_{k}^{1}{\big (}B; \vec{A}^{n-k}(\vec{\varphi}){\big )}                               \\ [0.5 ex]
\Delta_B(A)(x_{2}) & = & \sum_{k=1}^n (-1)^{k-1} \delta_{k}^{2}{\big (}B; \vec{A}^{n-k}(\vec{\varphi}){\big )}                               \\
\;\vdots           &   &                                                                                                                     \\
\Delta_B(A)(x_{n}) & = & \sum_{k=1}^n (-1)^{k-1} \delta_{k}^{n}{\big (}B; \vec{A}^{n-k}(\vec{\varphi}){\big )}.
\end{array}
\]
\end{theorem}

In the rest of this section we restrict our attention to some properties of the matrix $\hat{\mathstrut B}$.

The characteristic polynomial of the matrix $\hat{\mathstrut B}$ is
$$
\begin{array}{lcl}
\Delta_{\hat{\mathstrut B}}(\lambda)
& = &
det(\lambda I - \hat{\mathstrut B}) =
\left|
\begin{array}{rrrr}
\lambda-b_1  &        -b_1  & \ldots &        -b_1  \\
       -b_2  & \lambda-b_2  & \ldots &        -b_2  \\
\vdots\;\;\; & \vdots\;\;\; & \ddots & \vdots\;\;\; \\
       -b_n  &        -b_n  & \ldots & \lambda-b_n
\end{array}
\right|                                                    \\[5ex]
& = &
\left|
\begin{array}{cccc}
\lambda-\sum_{i=1}^n b_i  & \lambda-\sum_{i=1}^n b_i & \ldots & \lambda-\sum_{i=1}^n b_i  \\
       -b_2               & \lambda-b_2              & \ldots &        -b_2               \\
\vdots                    & \vdots                   & \ddots & \vdots                    \\
       -b_n               &        -b_n              & \ldots & \lambda-b_n
\end{array}
\right|
 =
(\lambda-\sum_{i=1}^n b_i)
\left|
\begin{array}{cccc}
1       & 1           & \ldots & 1           \\
-b_2    & \lambda-b_2 & \ldots & -b_2        \\
\vdots  & \vdots      & \ddots & \vdots      \\
-b_n    & -b_n        & \ldots & \lambda-b_n
\end{array}
\right|                                                   \\[5ex]
& = &
(\lambda-\sum_{i=1}^n b_i)
\left|
\begin{array}{cccc}
1       & 1       & \ldots & 1        \\
0       & \lambda & \ldots & 0        \\
\vdots  & \vdots  & \ddots & \vdots   \\
0       & 0       & \ldots & \lambda
\end{array}
\right|
=\lambda^{n-1}(\lambda-\sum_{i=1}^n b_i).
\end{array}
$$

Since
$$\hat{\mathstrut B} \cdot \left(\hat{\mathstrut B} - \sum_{i=1}^n b_i I\right) =
\left[
\begin{array}{cccc}
b_1    & b_1    & \ldots & b_1    \\
b_2    & b_2    & \ldots & b_2    \\
\vdots & \vdots & \ddots & \vdots \\
b_n    & b_n    & \ldots & b_n
\end{array}
\right]
\cdot
\left[
\begin{array}{cccc}
-\sum_{i=2}^n b_i & b_1                        & \ldots & b_1    \\
b_2               & -\sum_{\hspace{-0.18cm}
                           \scriptsize
                          {\begin{array}{c}
                           i=1       \\[-0.3ex]
                           i \neq 2
                           \end{array}}}^n
                           \hspace{-0.2cm}b_i  & \ldots & b_2    \\
\vdots             & \vdots                    & \ddots & \vdots \\
b_n                & b_n                       & \ldots & -\sum_{i=1}^{n-1} b_i
\end{array}
\right]
=\mathbb{O},
$$
the minimal polynomial of the matrix $\hat{\mathstrut B}$ is $\mu_{\hat{\mathstrut B}}(\lambda)=\lambda(\lambda-\sum_{i=1}^n b_i)$.

\bigskip
The further consideration is divided into two parts.

\bigskip
\textbf{Part 1:} Let $\sum_{i=1}^n b_i \neq 0$.
Then the invariant factors of the matrix $\hat{\mathstrut B}$ are $\underbrace{\lambda, \ldots \lambda}_{n-2}, \lambda^2-\sum_{i=1}^n b_i\lambda$
and its elementary divisors are $\underbrace{\lambda, \ldots \lambda}_{n-1}, \lambda-\sum_{i=1}^n b_i$.
The geometric multiplicity of an eigenvalue $\lambda = 0$ is equal to its algebraic multiplicity.
So the Jordan and the rational canonical form of the matrix $\hat{\mathstrut B}$ are
$$
J=
\left[
\begin{array}{ccccc}
0      & 0      & \ldots & 0      & 0      \\
0      & 0      & \ldots & 0      & 0      \\
\vdots & \vdots & \ddots & \vdots & \vdots \\
0      & 0      & \ldots & 0      & 0      \\
0      & 0      & \ldots & 0      & \sum_{i=1}^n b_i
\end{array}
\right]
=diag(0, 0, \ldots, 0, \sum_{i=1}^n b_i)
\;\;\mbox{and}\;\;
C=
\left[
\begin{array}{ccccc}
0      & 0      & \ldots & 0      & 0      \\
0      & 0      & \ldots & 0      & 0      \\
\vdots & \vdots & \ddots & \vdots & \vdots \\
0      & 0      & \ldots & 0      & 1      \\
0      & 0      & \ldots & 0      & \sum_{i=1}^n b_i
\end{array}
\right].
$$
The eigenvectors corresponding to the eigenvalue $\lambda=0$ are nonzero solutions of the homogeneous system
$$
\begin{array}{lcl}
b_1x_1 + b_1x_2 + \ldots + b_1x_n   & = & 0      \\[0.5 ex]
b_2x_1 + b_2x_2 + \ldots + b_2x_n   & = & 0      \\
\;\vdots   &   &                                 \\
b_nx_1 + b_nx_2 + \ldots + b_nx_n   & = & 0.
\end{array}
$$
Therefore, we can take
$$v_1=\left[
\begin{array}{r}
 -1    \\
  1    \\
  0    \\
\vdots \\
  0
\end{array}
\right],\;
v_2=\left[
\begin{array}{r}
 -1    \\
  0    \\
  1    \\
\vdots \\
  0
\end{array}
\right],\;
\ldots \;
,\;
v_{n-1}=\left[
\begin{array}{r}
 -1    \\
  0    \\
  0    \\
\vdots \\
  1
\end{array}
\right]
$$
to be basis of the eigenspace of the matrix $\hat{\mathstrut B}$ that corresponds to the eigenvalue $\lambda=0$.
The eigenvector corresponding to the eigenvalue $\lambda=\sum_{i=1}^n b_i$ is nontrivial solution of the homogeneous system
$$
\begin{array}{rcl}
-\sum_{i=2}^n b_i x_1 + b_1 x_2 + \ldots + b_1 x_n & = & 0                       \\
              b_2 x_1 -\sum_{\hspace{-0.18cm}
                           \scriptsize
                          {\begin{array}{c}
                           i=1       \\[-0.3ex]
                           i \neq 2
                           \end{array}}}^n
                           \hspace{-0.2cm}b_i x_2 + \ldots + b_2 x_n & = & 0    \\
\;\vdots                                                             &   &      \\
              b_n x_1 + b_n x_2 + \ldots -\sum_{i=1}^{n-1} b_i x_n   & = & 0.
\end{array}
$$
Multiplying the first equation successively with $-b_2, -b_3, \ldots, -b_n$ and adding to the second, the third, and finally the last equation
multiplied by $b_1$, we obtain the equivalent system
$$
\begin{array}{rcl}
b_1 x_2 & = & b_2 x_1 \\
b_1 x_3 & = & b_3 x_1 \\
\;\vdots              \\
b_1 x_n & = & b_n x_1.
\end{array}
$$
Hence, we conclude that the eigenvector corresponding to the eigenvalue $\lambda=\sum_{i=1}^n b_i$ is
$v_n = [b_1 \; b_2 \; \ldots \; b_n]^T$.
Let us denote by $S$ the modal matrix formed by the eigenvectors $v_1, v_2, \ldots , v_n$.
The matrix $S$ is of the form
$$
\left[
\begin{array}{rrrrr}
-1     & -1     & \ldots & -1     & b_1       \\
 1     &  0     & \ldots &  0     & b_2       \\
\vdots & \vdots & \ddots & \vdots & \vdots    \\
 0     &  0     & \ldots &  0     & b_{n-1}   \\
 0     &  0     & \ldots &  1     & b_n
\end{array}
\right],
$$
i.e. it is doubly companion matrix.
Applying Lemma 2.1 and Theorem 3.1 from Wanicharpichat's paper \cite{Wanicharpichat2012},
we easy conclude that the determinant and the inverse of the matrix $S$ are
$$det(S)=(-1)^{n-1}\sum_{i=1}^n b_i$$
and
$$
S^{-1}=
-\frac{1}{\sum_{i=1}^n b_i}
\left[
\begin{array}{ccccc}
b_2    & -\sum_{\hspace{-0.18cm}
                           \scriptsize
                          {\begin{array}{c}
                           i=1       \\[-0.3ex]
                           i \neq 2
                           \end{array}}}^n
                           \hspace{-0.2cm}b_i  & \ldots & b_2 & b_2 \\
b_3    & b_3    & \ldots & b_3    & b_3                             \\
\vdots & \vdots & \ddots & \vdots & \vdots                          \\
b_n    & b_n    & \ldots & b_n    & -\sum_{i=1}^{n-1} b_i           \\
-1     & -1     & \ldots & -1     & -1
\end{array}
\right].
$$
For the matrices $S$ and $S^{-1}$ the equality $J=S^{-1}\hat{\mathstrut B}S$ holds.
Furthermore, it is easy to check that for the matrix
$$
R=
\left[
\begin{array}{ccccc}
1      & 0      & \ldots & 0      & 0                             \\
0      & 1      & \ldots & 0      & 0                             \\
\vdots & \vdots & \ddots & \vdots & \vdots                        \\
0      & 0      & \ldots & 1      & -\frac{1}{\sum_{i=1}^n b_i}   \\
0      & 0      & \ldots & 0      & 1
\end{array}
\right]
$$
the equality $C=R^{-1}JR$ holds.
Hence, for the matrix $T=SR$, we have $C=T^{-1}\hat{\mathstrut B}T$.

\bigskip
\textbf{Part 2:} We now consider the case $\sum_{i=1}^n b_i=0$.
The geometric multiplicity of the only eigenvalue $\lambda = 0$ of the matrix $\hat{\mathstrut B}$ is still equal to $n-1$,
but its algebraic multiplicity is $n$. In this case the Jordan and the rational canonical form are the same
$$
J=C=
\left[
\begin{array}{ccccc}
0      & 0      & \ldots & 0      & 0      \\
0      & 0      & \ldots & 0      & 0      \\
\vdots & \vdots & \ddots & \vdots & \vdots \\
0      & 0      & \ldots & 0      & 1      \\
0      & 0      & \ldots & 0      & 0
\end{array}
\right].
$$
Since $\hat{\mathstrut B}^2$ is equal to $\mathbb{O}$, we can take for a generalized eigenvector
column $u=[1 \; 0 \; \ldots \; 0]^T$.
Corresponding eigenvector in associate chain is
$v=\hat{\mathstrut B}u=[b_1 \; b_2 \; \ldots \; b_n]^T$.
In this case the transition matrix is
$$
T=
\left[
\begin{array}{rrrrrr}
-1      & -1     & \ldots & -1     & b_1     & 1       \\
 1      &  0     & \ldots &  0     & b_2     & 0       \\
 0      &  1     & \ldots &  0     & b_3     & 0       \\
\vdots  & \vdots & \ddots & \vdots & \vdots  & \vdots  \\
0       & 0      & \ldots &  1     & b_{n-1} & 0       \\
0       & 0      & \ldots &  0     & b_n     & 0
\end{array}
\right].
$$

\section{Main Results}

After consideration of the Jordan and the rational canonical form of the matrix $\hat{\mathstrut B}$,
we can turn back to examination of the reduction process of the system (\ref{system}).

\begin{theorem}
\label{theorem Jordan}
If $\sum_{i=1}^n b_i \neq 0$, then the system {\rm (\ref{system})} can be transform into the system
\begin{equation}
\begin{array}{rcl}
A(y_1)     & =      & \psi_1                          \\
A(y_2)     & =      & \psi_2                          \\
           & \vdots       &                           \\
A(y_{n-1}) & =      & \psi_{n-1}                      \\
A(y_n)     & =      & \sum_{i=1}^n b_i y_n + \psi_{n},
\end{array}
\end{equation}
where the columns $\vec{y}$ and $\vec{\psi}$ are determined by
$\vec{y}=S^{-1}\vec{x}$ and $\vec{\psi}=S^{-1}\vec{\varphi}$,
for a nonsingular matrix $S$ such that $J=S^{-1}\hat{\mathstrut B}S$.
\end{theorem}

The proof follows immediately by Lemma \ref{Jordan form} and the Jordan canonical form of the system matrix $\hat{\mathstrut B}$.

\begin{theorem}
\label{theorem rational}
The system {\rm (\ref{system})} can be transform into the system
\begin{equation}
\begin{array}{rcl}
A(z_1)                               & =      & \nu_1                                             \\
A(z_2)                               & =      & \nu_2                                             \\
                                     & \vdots &                                                   \\
(A^2   - \sum_{i=1}^n b_iA)(z_{n-1}) & =      & \nu_n + A(\nu_{n-1})-\sum_{i=1}^n b_i\nu_{n-1}    \\
A(z_n) -\sum_{i=1}^n b_i z_n         & =      & \nu_{n},
\end{array}
\end{equation}
where the columns $\vec{z}$ and $\vec{\nu}$ are determined by
$\vec{z}=T^{-1}\vec{x}$ and $\vec{\nu}=T^{-1}\vec{\varphi}$,
for a nonsingular matrix $T$ such that $C=T^{-1}\hat{\mathstrut B}T$.
\end{theorem}

Proof:
Applying Lemma \ref{rational form} the system (\ref{system}) is equivalent to the system
$$
\begin{array}{rcl}
A(z_1)     & =      & \nu_1                          \\
A(z_2)     & =      & \nu_2                          \\
           & \vdots      &                           \\
A(z_{n-1}) & =      & z_n + \nu_{n-1}                \\
A(z_n)     & =      & \sum_{i=1}^n b_i z_n + \nu_{n}.
\end{array}
$$
The reduced system is obtained by acting of the operator $A-\sum_{i=1}^n b_i$ on the penultimate equation
and by substituting the expression $A(z_n) -\sum_{i=1}^n b_i z_n$ appearing on the right-hand side of the equation with $\nu_{n}$.
\qed

Theorems \ref{theorem Jordan} and \ref{theorem rational} in the exactly same forms hold for the system (\ref{transpose system}).
Nonsingular matrices $G$ and $H$ such that $J=G^{-1}\check{\mathstrut B}G$ and $C=H^{-1}\check{\mathstrut B}H$
can be calculate in a similar manner as matrices $S$ and $T$ for which $J=S^{-1}\hat{\mathstrut B}S$ and $C=T^{-1}\hat{\mathstrut B}T$ hold.
Matrix $G$ is a doubly companion matrix of the form
$$
\left[
\begin{array}{rrrrr}
-\frac{b_2}{b_1} & -\frac{b_3}{b_1} & \ldots & -\frac{b_n}{b_1} & 1       \\
 1               &  0               & \ldots &  0               & 1       \\
\vdots           & \vdots           & \ddots & \vdots           & \vdots  \\
 0               &  0               & \ldots &  0               & 1       \\
 0               &  0               & \ldots &  1               & 1
\end{array}
\right],
$$
and for matrix $H$ we have $H=RG$, for $R$ from the section 3.

From now on we are concerned how general technique from \cite{MalesevicTodoricJovovicTelebakovic2012} can be used in this special case.
First of all, we calculate the coefficients of the adjugate matrix of the characteristic matrix of the matrix $\hat{\mathstrut B}$.
The characteristic matrix of the matrix $\hat{\mathstrut B}$ is matrix $\lambda I -\hat{\mathstrut B}$.
The adjugate matrix of the matrix $\lambda I - \hat{\mathstrut B}$ is matrix
$adj(\lambda I - \hat{\mathstrut B})=\lambda^{n-1}B_0 + \lambda^{n-2}B_1 + \ldots + \lambda B_{n-2} + B_{n-1}$.
The coefficients $B_0, B_1, \ldots , B_{n-1}$ satisfy following recurrence
$$B_0=I, \quad B_k = B_{k-1} \hat{\mathstrut B} + d_k I, \quad 1 \leq k < n,$$
where $d_k$ is a coefficient of the characteristic polynomial $\Delta_{\hat{\mathstrut B}}(\lambda)$ in front of $\lambda^{n-k}$, see \cite{Gantmacher00}.
Thus, in this particular case we have
$$B_0=I, \; B_1 = B_0 \hat{\mathstrut B} + d_1 I = \hat{\mathstrut B} - \sum_{i=1}^n b_i I = \left[
\begin{array}{cccc}
-\sum_{i=2}^n b_i & b_1                        & \ldots & b_1    \\
b_2               & -\sum_{\hspace{-0.18cm}
                           \scriptsize
                          {\begin{array}{c}
                           i=1       \\[-0.3ex]
                           i \neq 2
                           \end{array}}}^n
                           \hspace{-0.2cm}b_i  & \ldots & b_2    \\
\vdots             & \vdots                    & \ddots & \vdots \\
b_n                & b_n                       & \ldots & -\sum_{i=1}^{n-1} b_i
\end{array}
\right],
\; B_2 = \ldots = B_{n-1} = \mathbb{O}.$$
Therefore, the adjugate matrix $adj(\lambda I-\hat{\mathstrut B})$ is equal to $\lambda^{n-1}I + \lambda^{n-2}(\hat{\mathstrut B}- \sum_{i=1}^n b_i I)$.
Applying Theorem \ref{Theorem adjugate matrix}, the system (\ref{system}) can be transformed into totally reduced system
$$\left(\Delta_{\hat{\mathstrut B}}(\vec{A})\right)(\vec{x})=\vec{A}^{n-1}(\vec{\varphi})+\left(\hat{\mathstrut B}- \sum_{i=1}^n b_i I\right)\vec{A}^{n-2}(\vec{\varphi}).$$

We have thus proved following theorem.

\begin{theorem}
\label{TRS1}
Linear system of the operator equations {\rm (\ref{system})} implies the system,
which consists of the higher order operator equations as follows
\begin{equation}
\tag{$\hat{5}$}
\label{TRS system}
\begin{array}{lcl}
\Delta_{\hat{\mathstrut B}}(A)(x_{1})   & = & A^{n-1}(\varphi_1) + \sum_{i=2}^n (b_1 A^{n-2}(\varphi_i) - b_i A^{n-2}(\varphi_1))         \\[0.5 ex]
\Delta_{\hat{\mathstrut B}}(A)(x_{2})   & = & A^{n-1}(\varphi_2) + \sum_{\hspace{-0.18cm}
                                                                         \scriptsize
                                                                        {\begin{array}{c}
                                                                         i = 1         \\[-0.3ex]
                                                                         i \neq 2
                                                                        \end{array}}}^n
                                                                        \hspace{-0.2cm}(b_2 A^{n-2}(\varphi_i) - b_i A^{n-2}(\varphi_2))   \\
\;\vdots                     &   &                                                                                               \\
\Delta_{\hat{\mathstrut B}}(A)(x_{n})   & = & A^{n-1}(\varphi_n) + \sum_{i=1}^{n-1} (b_n A^{n-2}(\varphi_i) - b_i A^{n-2}(\varphi_n)).
\end{array}
\end{equation}
\end{theorem}

The adjugate matrix of the characteristic matrix of the matrix $\check{\mathstrut B}$ is equal to
$\lambda^{n-1}I +\lambda^{n-2}(\check{\mathstrut B}- \sum_{i=1}^n b_i I)$.
Hence, the corresponding theorem for totally reduced system of system (\ref{transpose system}) reads.

\begin{theorem}
\label{TRS2}
Linear system of the operator equations {\rm (\ref{transpose system})} implies the system,
which consists of the higher order operator equations as follows
\begin{equation}
\tag{$\check{5}$}
\label{TRS transpose system}
\begin{array}{lcl}
\Delta_{\check{\mathstrut B}}(A)(x_{1})   & = & A^{n-1}(\varphi_1) + \sum_{i=2}^n b_i(A^{n-2}(\varphi_i) - A^{n-2}(\varphi_1))         \\[0.5 ex]
\Delta_{\check{\mathstrut B}}(A)(x_{2})   & = & A^{n-1}(\varphi_2) + \sum_{\hspace{-0.18cm}
                                                                           \scriptsize
                                                                          {\begin{array}{c}
                                                                           i = 1         \\[-0.3ex]
                                                                           i \neq 2
                                                                           \end{array}}}^n
                                                                           \hspace{-0.2cm}b_i(A^{n-2}(\varphi_i) - A^{n-2}(\varphi_2))   \\
\;\vdots                                  &   &                                                                                          \\
\Delta_{\check{\mathstrut B}}(A)(x_{n})   & = & A^{n-1}(\varphi_n) + \sum_{i=1}^{n-1} b_i(A^{n-2}(\varphi_i) - A^{n-2}(\varphi_n)).
\end{array}
\end{equation}
\end{theorem}

At the end of this section let us briefly analyze one application of the total reduction in~the~case when $A$ is a differential
operator $\frac{d}{dt}$ on the vector space of real functions in one unknown~$t$.
We consider systems (\ref{system})~and~(\ref{transpose system}) with the initial conditions $x_i(t_0) = c_i$, $1 \leq i \leq n$.
Recursively substituting the initial conditions into the systems we obtain $n\!\,-\,\!1$ additional initial condition
for each equation of the totally reduced systems (\ref{TRS system})~and~(\ref{TRS transpose system}).
Then Theorem \ref{TRS1} and \ref{TRS2} can be used for finding the unique solutions of the corresponding Cauchy problems.

\end{document}